\documentclass{article}

\usepackage{amssymb}
\usepackage{amsmath}

\def\nn{\nonumber}

\def\phi{\varphi}

\def\f{\frac}
\def\d{\partial}

\def\B{\Big}

\def\bs{\backslash}

\newtheorem{lemma}{Lemma}[section]
\newtheorem{theorem}{Theorem}

\usepackage{graphics}

\def\eqref#1{(\ref{eq.#1})}

\def\defined#1{{\em #1}}

\def\optional#1{}

\def\proof{\par\noindent{\bf Proof.\ }}

\def\eop{\vskip 3mm }

\newtheorem{proposition}{Proposition}

\def\hu{{\hat u}}

\def\bz{{\bf z}}
\def\bs{{\backslash}}

\def\D{{\bf D}}

\begin{document}

\title{Explicit enumeration of triangulations with multiple boundaries}
\author{
Maxim Krikun \\
Institut Elie Cartan\\
Universite Henri Poincare\\ 
Nancy, France \\
krikun@iecn.u-nancy.fr
}

\maketitle

\begin{abstract}
We enumerate rooted triangulations of a sphere with multiple holes 
by the total number of edges and the length of each boundary component.
The proof relies on a combinatorial identity due to W.T.~Tutte.
\end{abstract}

\section{Introduction}\label{sec.intro}
\subsection{Definitions}
A planar map is a class of equivalence of embedded graphs $G\hookrightarrow S^2$
by the homeomorphisms of $S^2$.  
We note by $V(G)$, $E(G)$ and $F(G)$ the sets of vertices, edges and faces 
of the the map $G$, respectively.

A \defined{map with holes}, is a pair $(G,H)$, $H\subset F(G)$, 
such that no two faces $h, h' \in H$ share a common vertex,
and all vertices at the boundary of $h_i \in H$ are distinct
(i.e. the boundary of $h_i$ is a cycle with no self-intersections).
In the following we refer to the faces $h\in H$ as \defined{holes}.
A map is called a \defined{triangulation}, if every face of $F(G)\bs H$ 
has degree $3$. If $H=\emptyset$, such triangulation is called a 
\defined{complete triangulation}.
In the following we will consider rooted triangulations,
that is triangulations with one distinguished directed edge, called the root.
In addition to that, we assume that the holes of a triangulation 
are enumerated by integers $0,\ldots,r$ and that the root is always located 
at the boundary of the $0$-th hole.

\subsection{Main result}
In this paper we solve explicitly the recursive equations for
generating functions planar triangulations with arbitrary number of holes,
in terms of the total number of edges and the length of each boundary
component.

The class of triangulations we consider is the most wide possible --- 
the underlying graph may contain multiple edges and loops.
Although this class is sometimes thought of as "pathological",
it turns out that the presence of loops is a feature
which greatly simplifies the calculations involved 
(e.g. compared to~\cite{K03}).

Our main result is the following
\begin{theorem}\label{T1}
Let $C_r(n,\alpha_0; \alpha_1,\ldots, \alpha_r)$ be the number of rooted
triangulations with $(r+1)$ hole, with $\alpha_j$ edges on the boundary
of the $j$-th hole and $n$ edges in total. Then we have,
letting $m = \alpha_0 + \ldots + \alpha_r$,
\begin{equation}\label{eq.T1}
C_r(n,\alpha_0; \alpha_1,\ldots, \alpha_r) = 
\f{4^k(2m+3k-2)!!}{(k+1-r)!(2m+k)!!} 
\,\,\alpha_0 \prod_{j=0}^r {2\alpha_j\choose \alpha_j},
\end{equation}
if $n=2m+3k$, and 
\[C_r(n,\alpha_0; \alpha_1,\ldots, \alpha_r) = 0 \]
if $n-2m\neq 0\,\mbox{(mod $3$)}$.
\end{theorem}

\subsection{Related work}
The case $r=0$ corresponds to the problem of enumeration of planar 
near-triangulations, solved by Tutte in~\cite{T62a}
using the method of recursive decomposition.
The same method, applied to the problem of enumeration of 
triangulations on an orientable surface of genus $g$, 
leads in a natural way to enumeration of triangulations (or maps) 
with multiple holes.
We were unable to obtain any general result in the non-planar case,
but for completeness we provide the corresponding recurrent relations 
in Section~\ref{sec.hg}, as well as the generating functions for 
the triangulations of orientable a torus and double torus ($g=1$ and $g=2$). 

The decomposition method used in our study and the equations involved 
are not new. 
The similar ideas were applied by Bender and Canfield (\cite{BC86}),
and later by Arqu\'es and Gioretti (\cite{AG}), 
to the asymptotical enumeration of arbitrary rooted maps on surfaces.

Similar equations appear 
under the name of {\em loop}, or {\em Schwinger-Dyson equations}
in some models of two-dimensional quantum gravity.
Ambj{\o}rn et al.~studied the asymptotical number of triangulations 
(and some more general classes of maps) on the sphere and higher genera surfaces
with multiple holes (see Chapter 4 in~\cite{ADJ}).
We have found that the Proposition~\ref{P1} in section~\ref{sec.proof} 
below looks very similar to the formula $(4.95)$ in~\cite{ADJ}, 
which expresses the generating function of planar maps with multiple boundary 
components via the repeated application of 
the so-called \defined{loop insertion} operator.


A simplified version of loop insertion operation may be described as follows.
Given a complete rooted triangulation, one can cut it along the root edge,
and identify the obtained hole with two edges of an additional triangle.
This operation provides a bijection between the complete rooted triangulations 
with $n$ edges, and triangulations with $n+2$ edges and a single hole 
of length $1$.
Thus taking $C_0(n+2,1)$, we recover the formula
\[ C_0(n+2,1) = \f{ 2 \cdot 4^{k-1} (3k)!! }{ (k+1)! (k+2)!! }, 
\qquad n=2+3k,
\]
which gives, by duality, the number of \defined{almost trivalent maps} with $k$ vertices 
(sequence A002005 in~\cite{SLO}), computed by Mullin, Nemeth and Schellenberg in~\cite{MNS}.

\subsection{Plan of the paper}
This paper is organized as follows. In section~\ref{sec.rel} we describe
the recursive decomposition procedure for triangulations and derive 
equations on the corresponding generating functions,
then solve explicitly these equations for $r=0,1,2,3$.
In Section~\ref{sec.extr} we calculate explicitly the coefficients $C_r$ for $r=0,2$.
This allows to suggest the main formula of Theorem~\ref{T1},
which is then proved in section~\ref{sec.proof}.
The proof closely follows that of~\cite{T62}.

\section{Recurrent relations}\label{sec.rel}
\subsection{Planar triangulations with holes}
Let $C_k(n,m;\alpha_1,\ldots,\alpha_k)$ be the number of rooted planar triangulations 
with $(k+1)$ holes $H=(h_0,h_1,\ldots,h_k)$, 
such that there are $m$ edges at the boundary of $h_0$, 
$\alpha_j$ edges on the boundary of $h_j$, $j=1,\ldots,k$ 
and $n$ edges total.

First, let us remind the recursive decomposition method.
Given a rooted planar triangulation $G$ with one hole
(that is, a triangulation of a disk),
and assuming that there is at least one triangle,
one deletes the triangle $t_0$ that contains the root.
In function of the position of a vertex $v$, opposite to the root edge in $t_0$,
there are two possibilities:
\begin{itemize}
\item[(A)] if $v$ is an internal vertex of the triangulation, one obtains a new triangulation
with one face less and one more edge on the boundary.
\item[(B)] if $v$ lies on the boundary of $G$, one cuts the resulting map in two parts,
with $(n_1,n_2)$ edges and the boundaries of length $(m_1,m_2)$,
such that $n_1 + n_2 = n - 1$ and $m_1 + m_2 = m + 1$,
$(n,m)$ being the number of edges and boundary length of the original configuration.
\end{itemize}
As the final object one obtains a planar map, consisting of a single edge, which we treat 
as a triangulation with $0$ faces, $1$ edge and one hole with boundary length $2$.

Now if $G$ is a triangulation with multiple ($k+1$, say) holes, 
there exists a third possibility for $v_0$, namely
\begin{itemize}
\item[(C)] if $v_0$ is located at the boundary of the hole $h_j$,
then after erasing the root edge one cuts the resulting map along $v_0$, 
obtaining a map with one hole less, and with boundary of $h_0$ having length
$m' = m + \alpha_j + 1$, $\alpha_j$ being the length of the boundary of $h_j$
in the original triangulation.
\end{itemize}
Now let $U_k(x,y,z_1,\ldots,z_k)$ be the multivariate generating function 
\[
U_k(x,y,z_1,\ldots,z_k)
= \sum_{N\ge1} \sum_{m\ge1} \sum_{\alpha_j\ge1} C(n,m;\alpha_1,\ldots,\alpha_k) x^n y^m z_1^{\alpha_1} \cdots z_k^{\alpha_k}.
\]
Translating the above decomposition procedure 
into the language of generating function,
we get the following 
\begin{lemma}
The following equations hold
\begin{eqnarray}
U_0(x,y) &=& x y^2 + \f{x}{y}\B( U_0(x,y) - y L_0(x)\B) + \f{x}{y} U_0^2(x,y) \label{eq.U0} \\
U_k(x,y;\bz) &=& 
\f{x}{y}\B( U_k(x,y;\bz) - y L_k(x;\bz) \B)       \nonumber\\
&&{} + \f{x}{y} \sum_{\omega\subset I_k} U_{|\omega|}(x,y;\bz_{\omega}) U_{k-|\omega|}(x,y;\bz_{I_k\bs\omega})      \nonumber \\
&&{} + \sum_{j=1}^k \B[
\f{x}{y-z_j} \B( \f{z_j}{y} U_{k-1}(x,y;\hat\bz_j) - \f{y}{z_j} U_{k-1}(x,z_j;\hat\bz_j) \B)
\nonumber \\
&& \phantom{{} + \sum_{j=1}^k\B[}
   + x L_{k-1}(x;\hat\bz_j) \B]
\label{eq.Uk}
\end{eqnarray}
where
\[ L_k(x;\bz) = [y] U_k(x,y;\bz), \]
$I_k = \{1,2,\ldots,k\}$ is the index set, 
the sum is over all subsets $\omega$ of $I_k$ 
(including empty set and $I_k$ itself),
$\bz$ stands for $z_1,\ldots,z_k$, 
$\bz_{\omega}$ is the list of variables $z_j$ with $j\in\omega$,
and $\hat\bz_j$ stands for $z_1,\ldots,z_k$ without $z_j$.
\end{lemma}

\proof 
The equation~\eqref{U0} is a classical relation for the generating function
of near-triangulations: the term $xy^2$ accounts for the special single-edged map,
the term, linear in $U_0$, corresponds to the case~(A) above, and the quadratic term
to the case~(B).

In~\eqref{Uk}, the first term on the right-hand side is derived exactly the
same way as in~\eqref{U0};
the summation over~$\omega$ corresponds to the possible ways to distribute
the $k$ enumerated holes between the two parts of a triangulation in case~(B).

To see how the summation over $j$ in~\eqref{Uk} arises, consider first the case $k=1$,
i.e. a triangulation with two holes. 
When the rule~(C) above applies, removing of the root edge merges the two holes,
of lengths $\alpha_0$ and $\alpha_1$, into a single hole of length $(\alpha_0+\alpha_1+1)$.
This gives the following contribution to $U_1(x,y,z)$:
\begin{eqnarray*}
\lefteqn{\sum_{n\ge0}\sum_{\alpha_0\ge1, \alpha_1\ge1} 
C_0(n-1,\alpha_0+\alpha_1+1) x^{n} y^{\alpha_0} z^{\alpha_1}} &&
\\
&=& x \sum_{n\ge 0} \sum_{m\ge3} C_0(n-1,m) x^{n-1} ( y z^{m-2} + y^2 z^{m-3} \ldots + y^{m-2} z)
\\
&=& x \sum_{n\ge 0} \sum_{m\ge3} C_0(n-1,m) x^{n-1} \f{ z y^{m-1} - y z^{m-1} }{y-z} \\
&=& \f{x}{y-z}\B[
    \f{z}{y} \B( U_0(x,y) - U_0(x,0) - y [t]U_0(x,t) - y^2 [t^2] U_0(x,y)\B)
\\
&& \phantom{Cix}
   -\f{y}{z} \B( U_0(x,z) - U_0(x,0) - z [t]U_0(x,t) - z^2 [t^2] U_0(x,y)\B)\B]
\\
&=& \f{x}{y-z}\B( \f{z}{y} U_0(x,y) - \f{y}{z} U_0(x,z)\B) + x[t]U_0(x,t).
\end{eqnarray*}
A general case $k\ge1$ follows similarly, since when merging the hole $h_j$ with the 
hole $h_0$, all other holes remain intact.
\eop

\subsection{Solution of recurrent equations}
The equations \eqref{U0}, \eqref{Uk} may be solved exactly.
First, \eqref{U0} is solved using the quadratic method, giving
\begin{equation}
U_0(x,y) = \f{h-y}{2h}\sqrt{1-4h^2y} - \f{x-y}{2x},
\end{equation}
where $h = h(x)$ is a positive power series in $x$, 
satisfying the relation
\begin{equation}\label{eq.h}
8 h^3 x^2-h^2+x^2 = 0,
\end{equation}
namely
\begin{equation}\label{eq.hseries}
h(x)
= \sum_{k=0}^\infty \f{4^k (3k-1)!!}{k! (k+1)!!} x^{3k+1}
= x \sum_{k=0}^\infty \f{2^k (3k-1)!!}{k! (k+1)!!} (2x^3)^k
\end{equation}
(cf. sequence A078531 in~\cite{SLO}).

Next, one may solve \eqref{Uk} with respect to $L_k(x;\bz)$ and 
group the terms containing $U_k(x,y;\bz)$, obtaining 
\begin{equation}\label{eq.Ukf}
x L_k(x,t;\bz) = \f1y \B(x-y+2x U_0(x,y)\B) U_k(x,y;\bz) + W_k(x,y;\bz),
\end{equation}
where $W_k(x,y;\bz)$ is the sum of terms in \eqref{Uk}, not containing $U_k$,
\begin{eqnarray*}
W_k\lefteqn{(x,y;\bz) =
\f{x}{y} \sum_{\omega\subset I_k\atop 1 < |\omega| <k} 
         U_{|\omega|}(x,y;\bz_{\omega}) 
         U_{k-|\omega|}(x,y;\bz_{I_k\bs\omega}) 
	 }&&\\
&&{} + \sum_{j=1}^k \B[
\f{x}{y-z_j} \B( \f{z_j}{y} U_{k-1}(x,y;\hat\bz_j) - \f{y}{z_j} U_{k-1}(x,z_j;\hat\bz_j) \B)
   + x L_{k-1}(x;\hat\bz_j) \B].
\end{eqnarray*}
Note that the left hand side of \eqref{Ukf} does not depend on $y$,
and the factor $(x-y+2x U_0(x,y)) = (x/h)(h-y)\sqrt{1-4h^2y}$ 
vanishes at $y=h$, thus
\begin{equation}
U_k(x,y;\bz) = hy \f{ W_k(x,h;\bz) - W_k(x,y;\bz) }{ x (h-y)\sqrt{1-4h^2y} }.
\end{equation}
In particular, we have
\begin{eqnarray}
U_1(x,y;z) &=& 
\f12 \f{ z (1- \sqrt{1-4h^2y}) - y (1 - \sqrt{1-4h^2z}) } { (y-z) \sqrt{1-4h^2y}}, \\
U_2(x,y;z_1,z_2) &=& 
\f{ 8h^5 y (1-\sqrt{1-4h^2z_1}) (1-\sqrt{1-4h^2z_2}) }
  { (1-4h^3) (1-4h^2y)^{3/2} \sqrt{1-4h^2z_1}\sqrt{1-4h^2z_2}}
\end{eqnarray}
It is somewhat more convenient to consider the "symmetrized" functions
\begin{equation}
U_k^{sym}(x,y;z_1,\ldots,z_k) = z_1\cdots z_k \,\f{\d^k}{\d z_1 \cdots \d z_k} U_k(x,y;z_1,\ldots,z_k),
\end{equation}
which correspond to adding an additional root on each of the $k$ holes $h_1,\ldots,h_k$.
The functions $U_k^{sym}$ are then symmetric in $(y,z_1,\ldots,z_k)$:
\begin{equation}
U_1^{sym}(x,y;z) = 
\f{4h^4 yz}
{\B(\sqrt{1-4h^2y} + \sqrt{1-4h^2z}\B)^2 \sqrt{1-4h^2y} \sqrt{1-4h^2z} },
\end{equation}
\begin{equation}
U_2^{sym}(x,y;z_1,z_2) =
\f{ 32h^9 y z_1 z_2 }
  { (1-4h^3) (1-4h^2y)^{3/2} (1-4h^2z_1)^{3/2} (1-4h^2z_2)^{3/2}}
\end{equation}
\begin{equation}
U_3^{sym}(x,y;z_1,z_2,z_3) = 
\f{ 3072 h^{14} y z_1 z_2 z_3 \times P_3(h, y,z_1,z_2,z_3)} 
  { (1-4h^3)^3 (1-4h^2y)^{5/2} \displaystyle \prod_{j=1}^3 (1-4h^2z_j)^{5/2}},
\end{equation}
where 
\begin{eqnarray*}
P_3(h,y,z_1,z_2,z_3) &=& 
1 -  3\sigma_{(1)}(h^3, h^2 y, h^2 z_1, h^2 z_2, h^2 z_3) \\
&&{}  + 8\sigma_{(1,1)}(h^3, h^2 y, h^2 z_1, h^2 z_2, h^2 z_3) \\
&&{}  - 16\sigma_{(1,1,1)}(h^3, h^2 y, h^2 z_1, h^2 z_2, h^2 z_3) \\
&&{} + 256 h^{11} y z_1 z_2 z_3.
\end{eqnarray*}
and $\sigma_{(1)}$, $\sigma_{(1,1)}$, $\sigma_{(1,1,1)}$ are Schur 
polynomials.


\subsection{Triangulations of higher genera}\label{sec.hg}
The decomposition procedure extends naturally to the triangulations 
of genus $g$ with the following essential changes
(here, as above, $v_0$ denotes the vertex opposite to the rooted edge 
 in the triangle which is removed)
\begin{itemize}
\item[(C')]  If the vertex $v_0$ lies on the boundary of $h_0$, the map is separated into two parts,
       and both the holes and the genus should be distributed between these parts;
\item[(D)]  It is possible that $v_0$ lies on the boundary of $h_0$ in such 
       a way that after deleting the triangle $t_0$ and cutting the map 
       along $v_0$ the map stays connected 
       (imagine the hole $h_0$ wrapping around the torus). In such case
       the resulting map will have genus $g-1$ and one more hole.
\end{itemize}
Let $T_{g,k}(x,y;z_1,\ldots,z_k)$ be the generating function of triangulations of genus $g$ with $(k+1)$ hole (obviously, $T_{0,k} = U_k$).
The decomposition procedure leads to the recursive relations,
similar to the main equation in~\cite{BC86}.
\begin{lemma}
The following relations hold:
\begin{eqnarray}
T_{g,k}(x,y;\bz) &=& \f{x}{y}\B( T_{g,k}(x,y;\bz)  - y [t] T_{g,k}(x,t;\bz) \B)   \nn\\
&&{} + \f{x}{y} \sum_{i=1}^g
     \sum_{\omega\subset I_k} T_{i, |\omega|}(x,y;\bz_{\omega}) T_{g-i, k-|\omega|}(x,y;\bz_{I_k\bs\omega}) \nonumber \\
&&{} + \sum_{j=1}^k \B[
     \f{x}{y-z_j} \B( \f{z_j}{y} T_{g,k-1}(x,y;\hat\bz_j)
                    - \f{y}{z_j} T_{g,k-1}(x,z_j;\hat\bz_j) \B)
\nonumber\\		    
&& \phantom{{} + \sum_{j=1}^k\B[}
     + x [t] T_{g,k-1}(x,t;\hat\bz_j) \B]
     \nonumber\\
&&{} + x \f{\d}{\d t} T_{g-1,k+1}(x,y; z_1,\ldots, z_k, t) \B|_{t=y}.
\label{eq.Tk}
\end{eqnarray}
\end{lemma}
\proof 
When the case~(D) applies, after removing the root edge we get a triangulation
with an additional hole, and with a distinguished vertex on the boundary of 
this hole (the image of $v_0$).
This gives the last term in~\eqref{Tk}, and the rest is similar to~\eqref{Uk}.
\eop

The equation~\eqref{Tk} may be solved analogously to~\eqref{Uk}.
In particular, we find generating function
for triangulations of genus $1$ and $2$ with one hole
\begin{eqnarray}
T_{1,0}(x,y) &=& \f{ (1-16h^5y)h^5y} { (1-4h^3)^2 (1-4h^2y)^{5/2} } \\
T_{2,0}(x,y) &=& \f{ P_{2,0}(h,y) }{ (1-4h^3)^7 (1-4h^2y)^{11/2}},
\end{eqnarray}
where 
\begin{eqnarray*}
P_{2,0}(h,y) &=&
3h^{11}y (35+184h^3+48h^6 )  \\
&&\times ( 1024 h^{11} y^4 + 1024 h^{12} y^3 - 1280 y^3 h^9 + 1 ) \\
&&{}+ 128h^{18} y^3 (545 + 1488h^3 - 3216h^6 + 2560h^9) \\
&&{} + 64h^{16} y^2 ( -307-480h^{6}+256h^{9}+ 324h^{3} )
\end{eqnarray*}

\section{Extracting exact coefficients}\label{sec.extr}
\subsection{Lagrange inversion}
Letting $h = x\sqrt{1+\zeta}$ and $ t = x^3$ in \eqref{h} we get
\begin{equation}\label{eq.zt}
\zeta = 8t ( 1+\zeta)^{3/2},
\end{equation}
so the Lagrange's inversion theorem applies,
and we have, assuming $n=m+3k$,
\begin{eqnarray}
[x^n] h^m &=& [x^{n-m}] (h/x)^m = [t^k](1+\zeta)^{m/2} \nonumber\\
&=& 
\f{1}{k} [\lambda^{k-1}]\B\{ \f{m}{2}(1+\lambda)^{m/2-1} (1+\lambda)^{3k/2} \B\}
\nonumber\\
&=& \f{m}{k!} 4^k \f{(m+3k-2)!!}{(m+k)!!}.
\label{eq.hmn}
\end{eqnarray}
In particular this gives the formula~\eqref{hseries} for $h(x)$.

For $U_0$ we have the following series expansion in $y$
$$
U_0(x,y) =  \f{ h - x  + 2h^3 x}{2hx} y
+ \sum_{m=0}^\infty \f{1}{m+1}{2m\choose m}\B(1 - \f{4m+2}{m+2} h^3 \B) h^{2m+1} y^{m+2}.
$$
Letting $n=2m+3k$, $k\ge-1$ and using \eqref{hmn} we obtain
\begin{equation}
[x^n y^m] U_0(x,y) = 
m {2m\choose m} \f{4^k (2m+3k-2)!! } { (k+1)! (2m+k)!! },
\label{eq.u0coeff}
\end{equation}
and $[x^n y^m] U_0(x,y) = 0$ if $n+m\neq 0\,(mod\,3)$.

Now note that $U_2^{(sym)}$ has the product form, 
so the expansion is particularly easy to calculate.
First we'll need the coefficients
\begin{eqnarray*}
[x^n] \B\{ \f{32 h^9}{1-4h^3} (4h^2)^{m-3} \B\}
&=& \f12\sum_{j=0}^\infty 4^{m+j} [x^n] h^{2m+3j+3} \\
&=& \f18 2^{2m+2k}(2m+3k-2)!! \sum_{j=1}^k \f{2m+3j}{(k-j)!(2m+2j+k)!!} \\
&=& \f18 2^{2m+2k} \f{(2m+3k-2)!!}{(k-1)!(2m+k)!!}.
\end{eqnarray*}
where $n=2m+3k$. Then we obtain
\begin{eqnarray}
[x^n y^{\alpha_0} z_1^{\alpha_1} z_2^{\alpha_2}] 
     \lefteqn{ 
U_2^{(sym)}(x,y,z_1,z_2)  
     } && 
\nonumber \\
&=&
\prod_{i=0}^2 \f{(2\alpha_i-1)!!}{2^{\alpha_i-1} (\alpha_i-1)!} \cdot 
[t^n] \B\{ \f{32h^9}{1-4h^3} (4h^2)^{m-3} \B\} 
\nonumber \\
&=&
2^{m+2k}\f{(2m+3k-2)!!}{(k-1)!(2m+k)!!}
\f{ (2\alpha_0-1)!!  (2\alpha_1-1)!!  (2\alpha_2-1)!! }
  { (\alpha_0-1)!  (\alpha_1-1)! (\alpha_2-1)! },
\nonumber \\
&=&
\alpha_0 \alpha_1 \alpha_2
{2\alpha_0 \choose \alpha_0} 
{2\alpha_1 \choose \alpha_1} 
{2\alpha_2 \choose \alpha_2} 
\cdot
\f{2^{2k} (2m+3k-2)!!}{(k-1)!(2m+k)!!}
\label{eq.u2coeff}
\end{eqnarray}
where $m = \alpha_0 +\alpha_1 +\alpha_2$, $n = 2m+3k$;
the coefficient is $C_2(n,\ldots)$ is null if $n-2m\neq 0\,(mod\,3)$.

The formulae~\eqref{u0coeff}, \eqref{u2coeff} allow to conjecture the 
following general formula
\begin{equation}
[x^n z_0^{\alpha_0} z_1^{\alpha_1} \ldots z_k^{\alpha_2}] 
U_k^{(sym)}(x,z_0;\bz)  
 =
 \f{4^k(2m+3k-2)!!}{(k+1-r)!(2m+k)!!}
\prod_{j=0}^r \alpha_j {2\alpha_j\choose \alpha_j}
\label{eq.cr}
\end{equation}
where $m=\alpha_0 + \ldots + \alpha_r$ and $n=2m+3k$.

Clearly, this formula is equivalent to~\eqref{T1},
and it further agrees with the above expressions for $U_1^{(sym)}$ and $U_3^{(sym)}$
(as can be seen by calculating explicitly few first terms in the 
power series expansions of these functions).

\subsection{The combinatorial identity}
The above expression~\eqref{cr} resembles a 
formula obtained by Tutte in \cite{T62}, 
for the number of slicings with $k$ external faces of degrees 
$2n_1,\ldots,2n_k$
\begin{equation}\label{eq.T62}
\gamma(n_1,n_2,\ldots,n_k) = \f{(n-1)!}{(n-k+2)!}\prod_{i=1}^k \f{(2n_i)!}{n_i! (n_i-1)!}
\end{equation}
The proof of \eqref{T62} relies on the following 
combinatorial identity:
\begin{eqnarray}
\lefteqn{\sum_{\omega \subset I}
\D^{|\omega|-k} \{\lambda \cdot f_\omega\} \cdot
\D^{|\bar\omega|-l} \{\mu \cdot f_{\bar\omega}\}  
}&&\nonumber\\
&=& \sum_{\omega\subset I\atop |\omega|<k}
    \sum_{i=0}^{k-1-|\omega|}
    (-1)^i {|\bar\omega|-l \choose i}
      \D^{|\bar\omega|-l-i}\{ 
      \D^{-k+|\omega|+i} \{\lambda\cdot f_\omega\} \cdot \mu f_{\bar\omega} \}
\nonumber\\      
&& {} +
    \sum_{\omega\subset I\atop |\bar\omega|<l}
    \sum_{i=0}^{l-1-|\bar\omega|}
    (-1)^i {|\omega|-l \choose i}
      \D^{|\omega|-l-i}\{ 
      \lambda f_{\bar\omega} \cdot
      \D^{-l+|\bar\omega|+i} \{\mu f_{\bar\omega}\}\}.
\label{eq.ci}
\end{eqnarray}
where $I$ is the set $\{1,\ldots,r\}$;
$\lambda$, $\mu$, $f_1, \ldots f_r$ are arbitrary (sufficiently often differentiable)
functions of a single parameter, say $x$, 
$f_\omega$ denotes the product $f_\omega = \prod\{f_i | i\in\omega\}$,
and $\D$ stands for the differentiation in $x$.
Whenever $\D$ appears with negative index (which can only be 
in the left-hand side of~\eqref{ci}), it is to be treated as
an operation of repeated integration, and it is assumed that
the constants of integration are fixed in some way for every 
$X$ that appears as the argument to $\D^{-1}$, so that $\D^{-1}(X)$
is uniquely defined.

\section{Proof of~Theorem~1}\label{sec.proof}
{\em
The proof is organized as follows: first we interpret the formula~\eqref{cr}
in terms of generating functions $U_k^{(sym)}$.
Then we use the equation~\eqref{Uk} and the combinatorial identity~\eqref{ci}
to show by induction that all of the generating function have the required form.
}

Note that in ~\eqref{cr} 
\[ \f{4^k(2m+3k-2)!!}{(k+1-r)!(2m+k)!!}
= \f{k!}{(k+1-r)!} \f1{2m} [x^{2m+3k}] h^{2m},
\]
thus we have (with $n=2m+3k$)
\begin{eqnarray}
C_r^{(sym)}\lefteqn{(n,\alpha_0;\alpha_1,\ldots,\alpha_r)}&&\nonumber \\
&=& \f{4^k(2m+3k-2)!!}{(k+1-r)!(2m+k)!!}
    \prod_{j=0}^r \alpha_j {2\alpha_j\choose \alpha_j} \nonumber \\
&=& \f{k!}{(k+1-r)!} [x^n] \f{h^{2m}}{2m}
    \prod_{j=0}^r \alpha_j {2\alpha_j\choose \alpha_j} \nonumber \\
&=& \f{k!}{(k+1-r)!} [x^n z_0^{\alpha_0}\cdots z_r^{\alpha_r}]
    \int\limits_0^{h^2} \prod_{j=0}^r \f{ 2s z_j }{(1-4s z_j)^{3/2}}\,\f{ds}{2s}
\label{eq.crint1}
\end{eqnarray}
since $ m = \alpha_0 + \ldots + \alpha_r$, and 
\[ \sum_{\alpha=0}^\infty \alpha {2\alpha \choose \alpha} z^\alpha 
   = \f{2z}{(1-4z)^{3/2}}.
\]
On the other hand, from~\eqref{hmn} we have 
\[ [x^{2m+3k}] h^{2m} = [t^k] (1+\zeta(t))^m \]
(where $\zeta(t)$ is defined by~\eqref{zt}),
thus we may continue~\eqref{crint1} with
\begin{eqnarray}
\ldots 
&=& \f{k!}{(k+1-r)!} [t^k] 
    \B\{ [z_0^{\alpha_0}\cdots z_r^{\alpha_r}]
    \int\limits_0^{1+\zeta(t)} 
    \prod_{j=0}^r \f{ 2s z_j }{(1-4s z_j)^{3/2}}\,\f{ds}{2s}
    \B\} \nonumber \\
&=& [t^{k+1-r}] \B(\f{\d}{\d t}\B)^{r-1}
    \B\{ [z_0^{\alpha_0}\cdots z_r^{\alpha_r}]
    \int\limits_0^{1+\zeta(t)} 
    \prod_{j=0}^r \f{ 2s z_j }{(1-4s z_j)^{3/2}}\,\f{ds}{2s}
    \B\}. 
\end{eqnarray}
Finally,~\eqref{cr} is equivalent to (assuming $r\ge1$)
\begin{eqnarray}
U_r^{(sym)}(x;z_0,\ldots,z_r) &=&
u_r^{(sym)}(x^3; x^2 z_0, \ldots, x^2 z_r), \nonumber \\
u_r^{(sym)}(t;z_0,\ldots,z_r) &=& 
t^{r-1} \B(\f{\d}{\d t}\B)^{r-1}
    \int\limits_0^{1+\zeta(t)} 
    \prod_{j=0}^r \f{ 2s z_j }{(1-4s z_j)^{3/2}}\,\f{ds}{2s}.
\end{eqnarray}
In the non-symmetric case, a similar calculation gives
\begin{eqnarray*}
U_r(x,y ;z_1,\ldots,z_r) &=&
u_r(x^3, x^2y; x^2 z_1,\ldots,x^2 z_r), \nonumber \\
u_r(t,y;z_1,\ldots,z_r) &=& 
t^{r-1} \B(\f{\d}{\d t}\B)^{r-1}
    \int\limits_0^{1+\zeta(t)} 
    \f{y}{(1-4s y)^{3/2}}
    \prod_{j=1}^r \B(\f{1}{\sqrt{1-4s z_j}} - 1\B)\,ds.
\label{eq.urint}
\end{eqnarray*}
Now if we put $\hu_k = t^{1-r} u_k$,
the statement of the Theorem~\ref{T1} is equivalent to the following
\begin{proposition}~\label{P1}
Let
\begin{equation}\label{eq.hurt}
\hu_r(t,y;z_1,\ldots,z_r)
= t^{1-r} U_r(t^{\f13}, t^{-\f23}y; t^{-\f23}z_1, \ldots, t^{-\f23} z_r). 
\end{equation}
Then for all $r\ge1$
\begin{equation}\label{eq.huki}
\hu_r(t,y;\bz)
 = \B(\f{\d}{\d t}\B)^{r-1}
    \int\limits_0^{1+\zeta(t)} 
    \f{y}{(1-4s y)^{3/2}}
    \prod_{j=1}^r \B(\f{1}{\sqrt{1-4s z_j}} - 1\B)\,ds.
\end{equation}
\end{proposition}
\proof
First, applying the transformation~\eqref{hurt} to $U_0$, $U_1$ we find
\[ \hu_0(t,y)
   = \f12\B(t- \f{y}{\sqrt{1+\zeta}}\B)\sqrt{1-4(1+\zeta)y} - \f{t-y}{2},
\]
\[ \hu_1(t,y,z)
   = \f{y\sqrt{1-4(1+\zeta)z}}{2(y-z)\sqrt{1-4(1+\zeta)y}}
     - \f{1}{2\sqrt{1-4(1+\zeta)y}} - \f{z}{2(y-z)}.
\]     
It can be verified by explicit integration that $\hu_1$ satisfies~\eqref{huki}.

Next, for all $r\ge2$ \eqref{huki} is equivalent to
\begin{equation}\label{eq.huki2}
\hu_k(t,y;\bz)
 = \B(\f{\d}{\d t}\B)^{r-2}
   \B\{
    \f{y\,\zeta'(t)}{(1-4(1+\zeta) y)^{3/2}}
    \prod_{j=1}^r \B(\f{1}{\sqrt{1-4(1+\zeta) z_j}} - 1\B) \B\}.
\end{equation}
From~\eqref{zt} we have 
\[ \zeta' = \f{16(1+\zeta)^{5/2}}{2-\zeta }, \]
so
\begin{eqnarray*}
\hu_2(t,y;z_1,z_2) &=& \f{8t(1+\zeta)^{5/2}}{1-\zeta/2}
 \f{y}{(1-4(1+\zeta)y)^{3/2}}  \\
&&\times \B(\f1{\sqrt{1-4(1+\zeta)z_1}}-1\B)
 \B(\f1{\sqrt{1-4(1+\zeta)z_2}}-1\B)
\end{eqnarray*}
satisfies~\eqref{huki} as well.

Now suppose that~\eqref{huki} holds for $r=0,1,\ldots,k-1$ for some $k\ge3$,
and let us show that it holds as well for $r=k$.

The equation~\eqref{Uk} leads to the following equation on $\hu_k$:
\begin{eqnarray}
\hu_k(t,y;\bz) &=& 
\f{t}{y}\B( \hu_k(t,y;\bz) - y\,\hat{l}_k(x;\bz) \B)       \nonumber\\
&&{} + \f{1}{y} \sum_{\omega\subset I_k} \hu_{|\omega|}(t,y;\bz_{\omega}) \hu_{k-|\omega|}(t,y;\bz_{I_k\bs\omega})      \nonumber \\
&&{} + \sum_{j=1}^k \B[ \f{1}{y-z_j} \B( \f{z_j}{y} \hu_{k-1}(t,y;\hat\bz_j) - \f{y}{z_j} \hu_{k-1}(t,z_j;\hat\bz_j) \B)  \nonumber \\
&&{} + \hat{l}_{k-1}(t;\hat\bz_j)\B],
\label{eq.hutk}
\end{eqnarray}
with
\[ \hat{l}_k(t;\bz) =  [y]\hu_k(t,y;\bz).  \]
Rewrite~\eqref{hutk} as
\begin{eqnarray}
\B(y-t-2\hu_0(t,y) \B)\,
\lefteqn{
     \hu_k(t,y,\bz) 
     }&& 
\nonumber \\
&=& \sum_{\omega\subset I_k \atop 1<|\omega|<k} 
                \hu_{|\omega|}(t,y;\bz_{\omega}) \hu_{k-|\omega|}(t,y;\bz_{I_k\bs\omega})
\nonumber \\
&&{} + \sum_{j=1}^k 
    \f{y}{y-z_j} \B( \f{z_j}{y} \hu_{k-1}(t,y;\hat\bz_j)
    - \f{y}{z_j} \hu_{k-1}(t,z_j;\hat\bz_j) \B)  
\nonumber \\
&&{} + y \B(\sum_{j=1}^k \hat{l}_{k-1}(t;\hat\bz_j) - t\,\hat{l}_k(t;\bz)\B).
\label{eq.hutk2}
\end{eqnarray}
In order to apply the combinatorial identity~\eqref{ci} to sum over~$\omega$ in~\eqref{hutk2},
we need to introduce some new notation.
We put 
\begin{eqnarray*}
\lambda &=& \f{y\cdot \zeta'}{(1-4(1+\zeta)y)^{3/2}}
           = \f{16y (1+\zeta)^{5/2}} {(2-\zeta)(1-4(1+\zeta)y)^{3/2}}, \\ 
f_j &=& \f{1}{\sqrt{1-4(1+\zeta)z_j}} - 1,
\end{eqnarray*}
let $\D = \displaystyle\f{\d}{\d t}$ and fix the following integrals
\begin{eqnarray*}
\D^{-1} \lambda &=& \f{1}{2\sqrt{1-4(1+\zeta)y}}, \\
\D^{-2} \lambda &=& \f{(\zeta - 8(1+\zeta)y) \sqrt{1-4(1+\zeta)y}}{16(1+\zeta)^{3/2}}, \\
\D^{-1}(\lambda f_j) &=& \f{y\sqrt{1-4(1+\zeta)z_j}}{2(y-z_j)\sqrt{1-4(1+\zeta)y}}
 - \f{1}{2\sqrt{1-4(1+\zeta)y}} - \f{z_j}{2(y-z_j)}.
\end{eqnarray*}
With these conventions we have
\[
\hu_0(t,y) = \D^{-2}(\lambda) - \f{(t-y)}{2},
\qquad
\hu_1(t,y,z_j) = \D^{-1}(\lambda f_j),
\]
and we have supposed that, according to~\eqref{huki2},
\[ \hu_{r}(t,y;\bz_\omega) = \D^{r-2}(\lambda f_\omega) \]
for all $\omega \subset I_k$ such that $r=|\omega|$, $2\le r\le k-1$.

Now applying~\eqref{ci} we obtain
\begin{eqnarray}
\sum_{\omega\subset I_k \atop 1<|\omega|<k} 
\lefteqn{
      \hu_{|\omega|}(t,y;\bz_{\omega}) \hu_{k-|\omega|}(t,y;\bz_{I_k\bs\omega})
   + 2 \D^{-2}(\lambda) \D^{k-2}(\lambda f_{I_k})
}&&\nonumber \\
&=& \sum_{\omega\subset I_k} 
    \D^{|\omega|-2}(\lambda f_\omega)
    \D^{|\bar\omega|-2}(\lambda f_{\bar\omega})
\nonumber \\
&=& 2\, \D^{k-2}\{ \D^{-2}(\lambda) \cdot \lambda f_{I_k} \}
   -2(k-2) \D^{k-3} \{\D^{-1}(\lambda) \cdot \lambda f_{I_k}\} \nonumber \\
&&{}  
+ 2\sum_{j=1}^{k} \D^{k-3} \{\D^{-1}(\lambda f_j) \cdot \lambda f_{I_k\backslash j}\} 
    \nonumber \\
&=& 2\, \D^{k-3}\B\{
\D^{-2}(\lambda) \cdot \D(\lambda) f_{I_k}
+ 3 \D^{-1}(\lambda) \cdot \lambda f_{I_k} \nonumber \\
&& 
{}
+ \sum_{j=1}^{k} 
  \B( 
  \D^{-2}(\lambda) \D(f_j) - \D^{-1}(\lambda)f_j + \D^{-1}(\lambda f_j) \B)
   \cdot \lambda f_{I_k\backslash j} \B\}
\label{eq.Dkj}
\end{eqnarray}
where in the last equality we used the identities
\[
\D\{\D^{-2}(\lambda) \lambda f_{I_k} \} = 
\D^{-1}(\lambda) \cdot \lambda f_{I_k} + 
\D^{-2}(\lambda) \D(\lambda) f_{I_k} + \sum_{j=1}^k \D^{-2}(\lambda) \D(f_j) 
\cdot \lambda f_{I_k\backslash j}
\]
and
\[
k\,\D^{-1}(\lambda) \cdot \lambda f_{I_k} = \sum_{j=1}^k \D^{-1}(\lambda) f_j 
\cdot \lambda f_{I_k\backslash j}.
\]
On the other hand, we have
\begin{eqnarray*}
\lefteqn{
\f{y}{y-z_j}\B(\f{z_j}{y} \hu_{k-1}(t,y,\hat\bz_j) - \f{y}{z_j}\hu_{k-1}(t,z_j,\hat\bz_j)\B)
}&&\\
&=& \D^{k-3}\B\{
    \B(
    \f{z_j}{y-z_j} 
    - 
    \f{y}{y-z_j} \f{(1-4(1+\zeta)y)^{3/2}}{(1-4(1+\zeta)z_j)^{3/2}}
    \B)
    \cdot \lambda f_{I_k\backslash j} 
    \B\}
\end{eqnarray*}
and 
\[ y - t - 2\hu_0(t,y) = -2\,\D^{-2}(\lambda), \]
so we further rewrite~\eqref{hutk2} as
\begin{eqnarray}
-2\D^{-2}\lefteqn{
(\lambda)\cdot \hu_k(t,y;\bz) + 2\D^{-2}(\lambda)\cdot \D^{k-2}(\lambda f_{I_k})
}&&\nonumber\\ 
&=& \D^{k-3}\B\{
      2\B(\D^{-2}(\lambda) \cdot \D(\lambda) f_{I_k}
          + 3 \D^{-1}(\lambda) \cdot \lambda f_{I_k} \B)
\nonumber\\
&&{}
+ \sum_{j=1}^k \lambda f_{I_k\backslash j} \cdot
  \B(
  2\D^{-2}(\lambda)\cdot\D(f_j) - 2\D^{-1}(\lambda)f_j + 2\D^{-1}(\lambda f_j) 
\nonumber\\
&&{}
  + \f{z_j}{y-z_j} 
    - 
    \f{y}{y-z_j} \f{(1-4(1+\zeta)y)^{3/2}}{(1-4(1+\zeta)z_j)^{3/2}}
    \B)
    \B\}
\nonumber\\
&&{} + y \B(\sum_{j=1}^k \hat{l}_{k-1}(t;\hat\bz_j) - t\,\hat{l}_k(t;\bz)\B).
\label{eq.hutk3}
\end{eqnarray}
A straightforward calculation then shows that
\begin{eqnarray}
\D^{-2}(\lambda)\cdot\D(\lambda) + 3(\D^{-1}\lambda)\cdot\lambda 
&=& \f{96 (1+\zeta)^{5/2}}{(2-\zeta)^2}y
\label{eq.calc1}
\end{eqnarray}
\begin{eqnarray}
2\D^{-2}
\lefteqn{
        (\lambda) \cdot \D(f_j) 
	- 2\D^{-1}(\lambda)\cdot f_j
	+ 2\D{(-1}(\lambda f_j) 
}&& \nonumber\\
&& {}+
\B(\f{z_j}{y-z_j}
 - \f{y}{y-z_j}
 \f{(1-4(1+\zeta)y)^{3/2}}{(1-4(1+\zeta)z)^{3/2}}\B) 
\nonumber\\
&=& \f{16 (1+\zeta)^{5/2} (\zeta-2+8(1+\zeta)z)}{(2-\zeta)^2(1-4(1+\zeta)z)^{3/2}}\, y 
\label{eq.calc2}
\end{eqnarray}
It follows from~\eqref{calc1}, \eqref{calc2} 
that the right-hand side of~\eqref{hutk3} is a linear function of $y$.
On the other hand, the left-hand side of~\eqref{hutk3} 
turns to zero both at $y=0$, because 
$\hu_k(t,0;\bz)=0$ and $\D^{2-k}\lambda f_{I_k} |_{y=0} =0$,
and at $y=\zeta/(8(1+\zeta))$, because then $\D^{-2}(\lambda)=0$.

Thus both sides of~\eqref{hutk3} are identically zero,
which implies $\hu_k = \D^{k-2}(\lambda f_{I_k})$.
This finishes the proof of both the proposition and Theorem~\ref{T1}.
\eop

\bibliographystyle{hplain} 
\bibliography{tri}

\end{document}